\documentclass[conference]{IEEEtran}
\IEEEoverridecommandlockouts

\usepackage{hyperref}
\usepackage{cite}
\usepackage{amsmath,amssymb,amsfonts}
\usepackage{graphicx}
\usepackage{textcomp}
\usepackage{color}
\def\BibTeX{{\rm B\kern-.05em{\sc i\kern-.025em b}\kern-.08em
		T\kern-.1667em\lower.7ex\hbox{E}\kern-.125emX}}

\begin{document}	
	\title{Simultaneously Parameter Identification and Measurement-Noise Covariance estimation of a Proton Exchange Membrane Fuel Cell\\
	}
	\author{
		\IEEEauthorblockN{1\textsuperscript{st} Razieh Ghaderi}
		\IEEEauthorblockA{\textit{Department of Mechatronic, Faculty of Engineering} \\
			\textit{Arak University}\\
			Arak, 38156-8-8349, Iran \\
			r-ghaderi@msc.araku.ac.ir, mer.gh1992@gmail.com}
		\and
		\IEEEauthorblockN{2\textsuperscript{nd}Abolghasem Daeichian}
		\IEEEauthorblockA{\textit{Department of Electrical Engineering, Faculty of Engineering},\\
			\textit{Arak University},\\
			Arak, 38156-8-8349, Iran \\
			\textit{Institute of advanced technology},\\
			\textit{Arak University},\\
			Arak, 38156-8-8349, Iran\\
			a-daeichian@araku.ac.ir, a.daeichian@gmail.com}
	}
	\maketitle
	
	\begin{abstract}
		This paper proposes the online parameters identification of semi-empirical models of Proton Exchange Membrane Fuel Cell (PEMFC). The covariance of unknown measurement noise is also estimated simultaneously. The actual data are fed to Kalman (for linear-in-parameters models) or extended Kalman filter (for nonlinear ones) which have been adapted for parameter identification. These filters suffer from the fact that the noise of the measurements is unknown. In order to tackle this conundrum, the measurement-noise is simultaneously estimated, and the estimation is used in the filters. The ultimate consequence of estimating measurement-noise is error reduction which has been demonstrated by simulation results. 
	\end{abstract}
	
	
	\section{Introduction}
	\label{section.introduction}
	Fossil fuels are not only air pollutant but also a limit source of energy. So, fuel cells (FC) have drawn lots of attention and has been rapidly developed \cite{R1,R2}. FC is an electrochemical device that produces electricity by converting fuel and oxidant that the result of this chemical reaction is to produce waste heat and water \cite{R3}.  
	Proton Exchange Membrane Fuel Cell (PEMFC), Alkaline fuel cell (AFC), and Phosphoric acid fuel cell (PAFC) are some kind of FCs. 		
	PEMFC has been ranked as a good source of power due to having not only high efficiency but also low high-power density. So, applying PEMFC in many applications such as an electric vehicle is on rise \cite{R4}.
	Several models have been proposed for PEMFC in literature, including black-box and Gray-box models. Black-box approach is widely employed to model fuel cells. For instance, feedforward neural networks \cite{R5}, recurrent neural networks \cite{R6}, nonlinear autoregressive with exogenous input (NARX) \cite{R7}, and radial basis function artificial neural networks \cite{R8} have been presented. Also, fuzzy logic control (FLC) as an adaptive neuro-fuzzy inference system (ANFIS) has been utilized in \cite{R9,R10}. A significant disadvantage of black-box models is that do not consider physical relationships. 
	
	Gray-box (semi-empirical) models propose an empirical equation which has some unknown parameters. For example, Squadrito et al. \cite{R11} and Kim et al. \cite{R12} proposed an equation between voltage and current of fuel cells. A pros of these models is considering the mass transport phenomena. Optimization methods such as recursive least square (RLS) \cite{R13,R14}, Kalmn filter (KF) and unscented Kalman filter (UKF) \cite{R15} have been conducted to identify the parameters of semi-empirical models. It worth to note that the measurement covariance have not been mentioned in these articles.
	
	This paper investigates the parameter identification of Squadrito and Kim models by simultaneously estimating the measurement-noise covariance. To this end, a learning equation has been proposed which exploit estimation algorithm given in \cite{R16}. The Squadrito model leads to an equation which is linear-in-parameters while the Kim causes a nonlinear one. Consequently, KF and EKF have been applied to the former and the latter, respectively. The results have been indicated that providing estimation of noise covariance to identification algorithm leads to less mean square error.
	
	The remainder of this paper is organized as follows: The Squadrito and Kim model have been presented in section \ref{section.model}. The algorithms of KF and EKF, which suit to parameter identification problem have been represented in section \ref{section.identification}. Section \ref{section.noisecovariance} is devoted to estimating noise covariance. The results are discussed in section \ref{section.results}. Finally, the paper has been concluded in section \ref{section.conclusion}.
	\section{PEMFC Models}
	\label{section.model}
	Several models have been introduced for PEMFCs. In this paper, two models are interested, which lead to linear- and nonlinear-in-parameters equations of regressors for identification. 
	
	\subsection{Squadrito et al. model}
	The first one is a static model proposed by Squadrito et al. \cite{R11}. In this model, the output voltage of FC is a function of its current as stated in the following equation: 
	\begin{equation}\label{Eq.Squadrito.model}
	v_{FC}=V_{0}-b\log(i_{FC})-r i_{FC}+\alpha i^{k}_{FC}\log(1-\beta i_{FC})
	\end{equation}
	where $v_{FC}$, $i_{FC}$, $V_{0}$, and $r$ are the output voltage of FC, the current of FC, the open-circuit voltage of the fuel cell, and the ohmic resistance of the fuel cell, respectively. $\alpha$ and $k$ are fitting parameters. $\beta$ is the inverse of limiting current $i_{L}$. $k$ is usually a constant between $1$ and $4$.
	An advantage of this model is taking mass transport over potential into accounts. 
	
	The measured values are $v_{FC}$ and $i_{FC}$ while the set of $\{V_{0}, r, \alpha, k\}$ are unknown parameters. Define
	\begin{eqnarray}\label{}
	\theta^T &=& \left[\theta_i\right]=\left[V_{0}, b, r, \alpha\right]\nonumber\\
	X^T &=& \left[x_i\right]=\left[1, -\log(i_{FC}), -i_{FC}, i^{k}_{FC}\log(1-\beta i_{FC})\right] \nonumber
	\end{eqnarray}
	where $x_i$ and $\theta_i$ are called regressors and regression coefficient (parameter) vectors, respectively. The $^T$ denotes the transpose of a vector. Then, the output $y(t)=v_{FC}(t)$ given by Eq.\ref{Eq.Squadrito.model} depends linearly on the parameters $\theta_i$ as:
	\begin{equation}\label{Eq.LinearRegression}
	y=\sum_{i=1}^{4}\theta_i x_i=X^T\theta.
	\end{equation}
	The aim is to identify the parameters $\theta$ such that minimize the sum of squared error loss function considering data samples $\{v_{FC}(t),i_{FC}(t)\}$, which are disturbed by a white noise.
	
	\subsection{Kim et al. model}
	A semi-empirical model to fit the output voltage and current of the PEMFC has been developed by Kim et al. \cite{R12}. This model is validated in different operating conditions including different temperature, pressure, and gas mixture. This model is given as:
	\begin{equation}\label{Eq.kim.model}
	v_{FC}=V_{0}-b\log (i_{FC})-r i_{FC}-m \exp (n i_{FC})
	\end{equation}
	where $v_{FC}$, $i_{FC}$, $b$ and $r$ are the PEMFC voltage, current, the Tafel parameters for oxygen reduction, and the ohmic resistance of the PEMFC, respectively. $V_{0}=V_{r}+b\log(i_{0})$ where $V_r$ and $i_0$ are the reversible potential of the cell and the exchange current. Also, $m$ and $n$ account for the mass transport phenomena as a function of current density. 
	
	Given Eq. (\ref{Eq.kim.model}), the unknown parameters can be considered as:
	\begin{eqnarray}\label{}
	\theta^T &=& \left[\theta_i\right]=\left[V_{0}, b, r, m, n\right]\nonumber
	\end{eqnarray}
	however, the output is not linear-in-parameters and a linear regression equation such as Eq.\ref{Eq.LinearRegression} can not be obtained. In this case, the Extended Kalman Filter has be employed for parameter identification.
	
	\section{Identification Algorithm}
	\label{section.identification}
	\subsection{Parameter identification of Squadrito model}
	There are many optimization algorithms which could be employed for online parameter identification such as the recursive least square (RLS). The Kalman filter, which is usually known as an optimal state estimator, could be applied as a parameter estimator. It extracts the interesting parameters from noisy observation. This filter predicts the next state variables and then updates the predicted value when the next measurement is received.
	In this section, the typical kalman filter and its extended version have been reviewed. Then, the parameter identification problem has been fitted to Kalman structure.

	Consider the true state and measurement model as:
	\begin{eqnarray}
	\theta(t+1) &=& F_t\theta(t)+B_tu(t)+w(t) \nonumber\\
	y(t) &=& H_t\theta(t)+n(t)
	\end{eqnarray}
	where $\theta$, $u$, $y$ are states, control-inputs, and measurements, respectively. The matrix $F_t$, $B_t$, and $H_t$ are state transition, control-input, and observations model. The process noise $w$ and the measurement noise $n$ are assumed to have zero-mean normal distribution with covariance $W_t$ and $R_t$, respectively.
	
	Typically, the Kalman filter has two alternate phase. The prediction phase as:
	\begin{eqnarray}
	\hat{\theta}(t|t-1) &=& F_t\hat{\theta}(t-1|t-1)+B_tu(t-1) \nonumber\\
	P(t|t-1) &=& F_tP(t-1|t-1)F_t^T+W_t
	\end{eqnarray}
	and the update phase as:
	\begin{eqnarray}
	\hat{\theta}(t|t) &=& \hat{\theta}(t|t-1)+K(t)\left(y(t)-H_t\hat{\theta}(t|t-1)\right) \nonumber\\
	P(t|t) &=& \left(I-K(t)H_t\right)P(t|t-1)\nonumber\\
	K(t) &=& P(t|t-1)H_t^T\left(H_tP(t|t-1)H_t^T+R_t\right)^{-1}
	\end{eqnarray}
	
	In the light of the fact that the parameters of PEMFC do not vary in time rapidly, and considering Squadrito model (Eq.\ref{Eq.LinearRegression}) as measurement model, the parameter estimation problem suits to state-space framework as: 
	\begin{eqnarray}\label{Eq.KF.model.statespace}
	\theta(t+1) &=& \theta(t)+w(t) \nonumber\\
	y(t) &=& X^T_t\theta(t)+n(t)
	\end{eqnarray}
	where $\theta(t)$ is the unknown parameters at discrete time $t$ which tend to be identified. $y(t)$ and $X_t$ are the measured  values and regressors, respectively. In addition, $w(t)$ and $n(t)$ are parameters and measurement noise with covariance matrices $W$ and $R$, respectively. The covariance matrices are typically chosen diagonal. The strength of parameters' time variance can be considered in the diagonal entry of $W$. This provides the capability to determine the forgetting factor for each parameter individually, which outperforms the RLS \cite{R20}. A forgetting factor equal to one is equivalent to $W=0$. Considering state and measurement equations (\ref{Eq.KF.model.statespace}), the Kalman filter algorithm becomes:
	\begin{eqnarray}\label{Eq.KF.algorithm}
	\hat{\theta}(t) &=& \hat{\theta}(t-1)+K(t)\left( y(t)-X^T_t\hat{\theta}(t-1)\right)\nonumber\\
	K(t) &=& \left(X^T_tP(t-1)X_t+R_t\right)^{-1}P(t-1)X_t\nonumber\\
	P(t) &=& \left(I-K(t)X^T_t\right)\left(P(t-1)+W\right)
	\end{eqnarray}                        
	where $P(t)$, $K(t)$, and $I$ are the error covariance matrix, the Kalman gain, and the identity matrix. Apparently, this filter requires noise covariance $R$. It commonly is considered as a constant.
	
	\subsection{Parameter identification of Kim model}
	The traditional Kalman filter has been extended for nonlinear state-space and measurement models by several approaches. Extended Kalman filter (EKF) is one of them linearizing the state-space model at each time step by Taylor series expansion. If the state and measurements evolve according to nonlinear but differentiable functions as:
	\begin{eqnarray}
	\theta(t+1) &=& f\left(\theta(t),u(t)\right)+w(t) \nonumber\\
	y(t) &=& h\left(\theta(t)\right)+n(t)
	\end{eqnarray}
	then the Extended Kalman filter prediction equation is:
	\begin{eqnarray}
	\hat{\theta}(t|t-1) &=& f\left(\hat{\theta}(t-1|t-1),u(t)\right) \nonumber\\
	P(t|t-1) &=& F_tP(t-1|t-1)F_t^T+W_t
	\end{eqnarray}
	The update phase is:
	\begin{eqnarray}
	\hat{\theta}(t|t) &=& \hat{\theta}(t|t-1)+K(t)\left(y(t)-h\left(\hat{\theta}(t|t-1)\right)\right) \nonumber\\
	P(t|t) &=& \left(I-K(t)H_t\right)P(t|t-1)\nonumber\\
	K(t) &=& P(t|t-1)H_t^T\left(H_tP(t|t-1)H_t^T+R_t\right)^{-1}
	\end{eqnarray}
	where
	\begin{eqnarray} 
	F_t&=&\frac{\partial f}{\partial \theta}\bigg |_{\hat\theta_{k-1|k-1},u_k}\nonumber\\
	H_t&=&\frac{\partial h}{\partial \theta}\bigg |_{\hat\theta_{k|k-1}}.
	\end{eqnarray}
		
	The Kim model has a nonlinear-in-parameters regression equation. Thus, EKF is a handy choice for its parameter identification. In this situation, the problem adapts to the state-space model: 
	\begin{eqnarray}\label{Eq.EKF.model.statespace}
	\theta(t+1) &=& \theta(t)+w(t) \nonumber\\
	y(t) &=& h\left(\theta(t)\right)+n(t)
	\end{eqnarray}
	where $h$ is a nonlinear function. The EKF algorithm is same as KF (Eq.\ref{Eq.KF.algorithm}), but $X(t)$ is the Jacobian of the function $h$ as:

	\begin{eqnarray}
	X_t &=& \frac{\partial h}{\partial \theta}\bigg |_{\hat\theta}
	\end{eqnarray}

\section{PEMFC Parameter Identification with Measurement Covariance Estimation}
\label{section.noisecovariance}
As stated in previous section, the KF and EKF can be utilized for parameter identification of the Squadrito and Kim models, respectively. However, the point is that the noise covariance $R_t$ should be determined. To this end, a straightforward and ordinary solution is to consider a constant value for $R_t$. Yet, more intellectual and rational approach is to provide an estimate for noise covariance in each time step. We propose a learning equation based on the noise covariance estimation given in \cite{R21}.

The covariance of the output-noise can be estimated as \cite{R22}, \cite{R23}:
\begin{equation}\label{Eq.Estimating.R}
   \hat{R}_t=\hat{l_{0}}-X_t(\hat{P}\hat X_t^{{T}})
\end{equation}
where $\hat{l_{0}}$ is an estimation of the innovation process $e_t=v_{FC}-\hat{v}_{FC}$ where $\hat{v}_{FC}$ is the estimation given by KF or EKF. So, the following learning equation has been proposed for $R_t$:
\begin{equation}
R_{t+1}=\lambda R_t+(1-\lambda)\hat R
\end{equation}
where $\lambda$ is a learning factor between $0$ and $1$. This equation learns the covariance of the measurement-noise from previous values and the current estimated value. The estimated value by Eq.\ref{Eq.Estimating.R} may diverge due to bad initial conditions. But, the learning algorithm considers the previous values of the measurement-noise covariance matrix and it became more hard to diverge. Also, a boundary could be imposed on the norm of measurement-noise covariance matrix to prevent divergence.

\section{RESULTS AND DISCUSSION}
\label{section.results}
\begin{figure}[!tbh]\
	\includegraphics[width=0.45\textwidth]{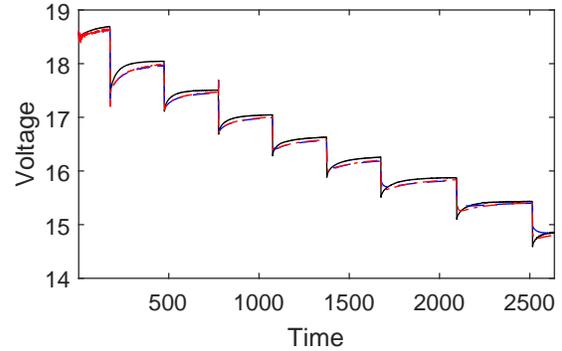}
	\caption{FC voltage: Solid line (black): experimental data; Dashed line (blue): KF with $R=1$; Dash-dot line (red): KF with estimating $R$}
	\label{Fig.Squadrito.voltage}
\end{figure}

A set of actual data consists of 26103 points with sampling interval as $0.1010$ that including $v_{FC}$ and $i_{FC}$ measured from a PEMFC with $36$ number of cells and $500$W rated power. KF is employed to identify the parameters of the Squdrito model. Two situations have been fulfilled: first, $R$ considered as a constant equal to $1$. Second, measurement noise has been estimated by the proposed method. The mean square error (MSE) of the predicted output voltage of the FC has been calculated for both situations, not only for all data but also for data after damping initial transient time. The results which have been shown in Fig.\ref{Fig.Squadrito.voltage} indicates that parameter identification with estimating $R$ outperforms considering constant $R$. 

\begin{figure}[!tbh]
	\includegraphics[width=0.45\textwidth]{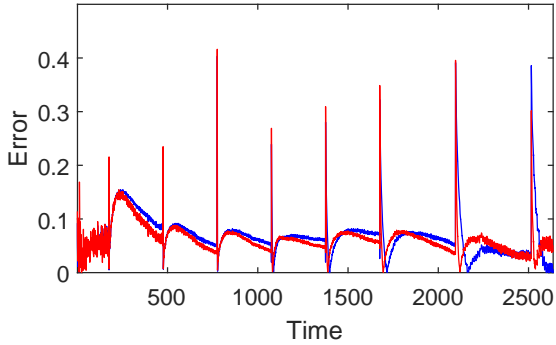}
	\caption{FC voltage estimation error}
	\label{Fig.Squadrito.error}
\end{figure}

\begin{figure}[!tbh]
	\includegraphics[width=0.4\textwidth]{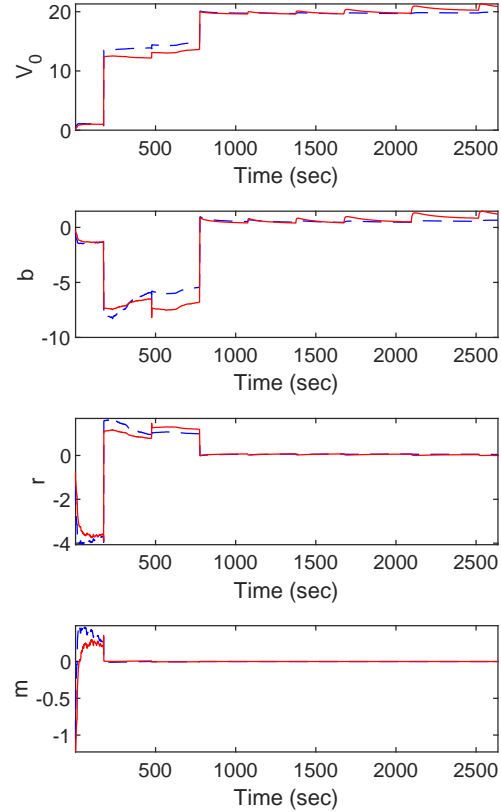}
	\caption{Estimated parameters of Squadrito model. Dashed line (blue): KF with $R=1$. Solid line (red): KF with estimating $R$}
	\label{Fig.Squadrito.parameters}
\end{figure}

The estimation error has been depicted in Fig.\ref{Fig.Squadrito.error}. Finally, The MSE of estimations in two situations has been given in Table \ref{Table.MSE}. Also, the estimated parameters are depicted in Fig. \ref{Fig.Squadrito.parameters}.


\begin{figure}[!tbh]\
	\includegraphics[width=0.45\textwidth]{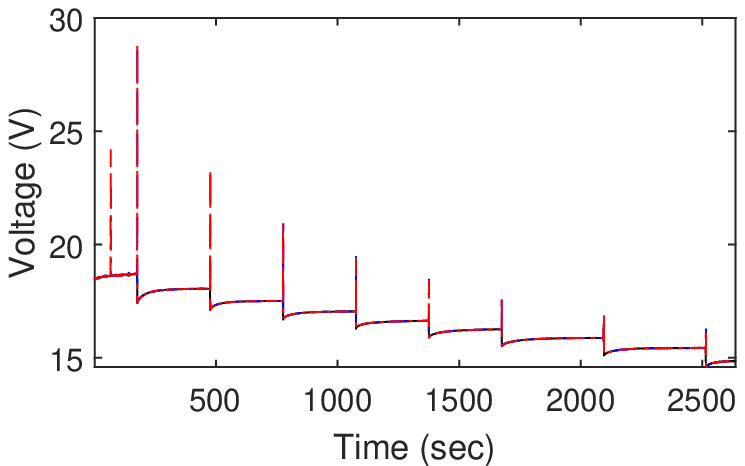}
	\caption{FC voltage: Solid line (black): experimental data; Dashed line (blue): EKF with $R=1$; Dash-dot line (red): EKF with estimating $R$}
	\label{Fig.Kim.voltage}
\end{figure}

\begin{figure}[!tbh]
	\includegraphics[width=0.45\textwidth]{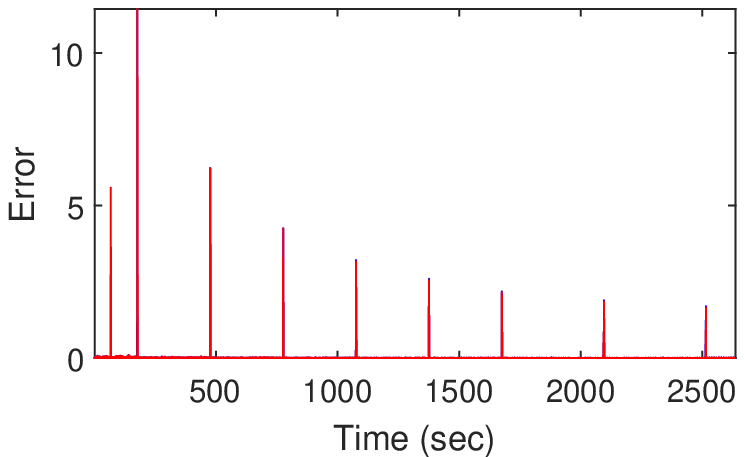}
	\caption{FC voltage estimation error}
	\label{Fig.Kim.error}
\end{figure}

\begin{figure}[!tbh]
	\includegraphics[width=0.4\textwidth]{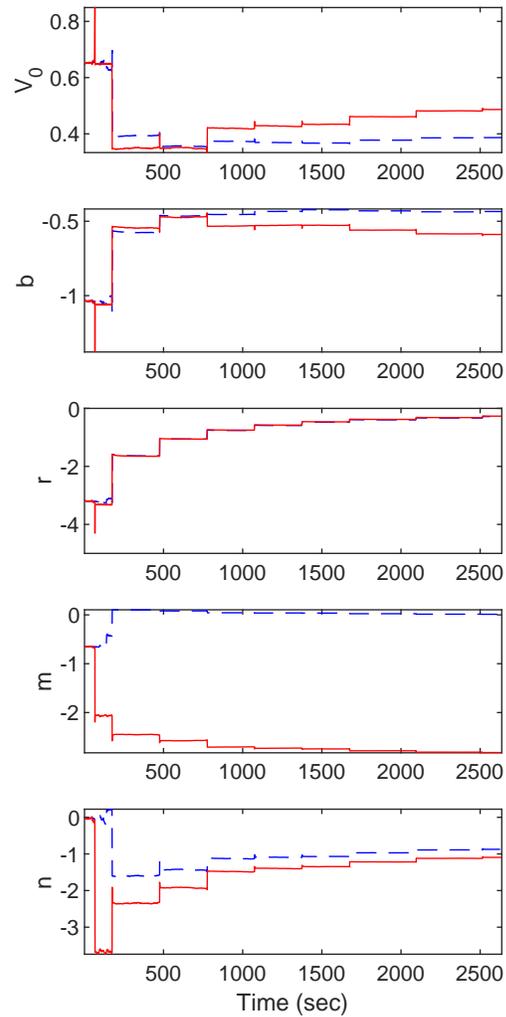}
	\caption{Estimated parameters of Kim model. Dashed line (blue): EKF with $R=1$. Solid line (red): EKF with estimating $R$}
	\label{Fig.Kim.parameters}
\end{figure}

Kim model is nonlinear-in-parameters. So, EKF is employed as an identification method. The actual and estimated output voltage of the FC are shown in Fig. \ref{Fig.Kim.voltage}. Also, the estimation error in two situations is displayed in Fig. \ref{Fig.Kim.error}. In this case, we have some peaks in estimation error when a sudden change of current occurs. The estimated parameters are depicted in Fig. \ref{Fig.Kim.parameters}.


\begin{table*}
	\begin{center}
		\caption{MSE of identification\label{Table.MSE}}
		\begin{tabular}{ccccccc}
			\hline\hline\\
			& & & \multicolumn{2}{c}{Constant R} & \multicolumn{2}{c}{Estimating R}\\\vspace{0.1cm} 
			Model & Algorithm & Parameters & MSE (1) & MSE (2) & MSE (1) & MSE (2)\\
			\hline\hline\\\vspace{0.1cm} 
			Squadrito & KF  & $\left[V_{0}, b, r, \alpha\right]$ & $6.176\times 10^{-3}$ & $5.404\times 10^{-3}$ & $4.643\times 10^{-3}$ & $3.866\times 10^{-3}$\\
			Kim       & EKF & $\left[V_{0}, b, r, m, n\right]$ & $2.403\times 10^{-2}$ & $4.484\times 10^{-3}$ & $2.393\times 10^{-2}$ & $4.431\times 10^{-3}$\\
			\hline
		\end{tabular}
	\end{center}
\end{table*}
By the way, the results demonstrate that estimating the covariance of measurement-noise decrease the MSE of parameter estimation for the PEMFC models.

\section{conclusion}
\label{section.conclusion}
The parameter identification of Squadrito and Kim models of PEMFC by simultaneously estimating the covariance of measurement-noise $R$ was investigated. A learning equation to estimate $R$ was proposed which try to learn from the estimation given in \cite{R22}. Also, Squadrito and Kim models result to equations which are linear and nonlinear in parameters. Thus, KF and EKF were selected for identification of Squadrito and Kim models, respectively. The results showed that estimating $R$ simultaneous with parameter identification outperform considering a constant $R$.

\end{document}